\documentclass[5p]{elsarticle}
\usepackage[makeroom]{cancel}
\usepackage{bm}
\usepackage{amsmath}
\usepackage{amsfonts}
\usepackage{siunitx}
\usepackage{hyperref}
\newcommand{\be}{\begin{equation}}
\newcommand{\ee}{\end{equation}}
\newcommand{\ber}{\begin{eqnarray}}
\newcommand{\eer}{\end{eqnarray}}
\newcommand{\p}{\partial}

\newcommand{\mat}[1]{\underline{\underline{#1}}}
\newcommand{\bbbr}{\mathbb{R}}

\begin{document}
\bibliographystyle{elsarticle-num}

\title{Large-Scale Magnetostatic Field Calculation in Finite Element Micromagnetics with $\mathcal{H}^2$-Matrices}
\author[rh]{Riccardo Hertel\corref{cor1}}
\address[rh]{Universit{\'e} de Strasbourg et CNRS, Institut de Physique et Chimie des Mat{\'e}riaux de Strasbourg, UMR 7504, 67000 Strasbourg, France}
\cortext[cor1]{Corresponding author}
\ead{riccardo.hertel@ipcms.unistra.fr}

\author[sb]{Sven Christophersen}

\author[sb]{Steffen B{\"o}rm}
\address[sb]{Department of Computer Science, 
Christian-Albrechts-Universit{\"a}t zu Kiel, 24118 Kiel, Germany}

\begin{frontmatter}
\begin{abstract}
Magnetostatic field calculations in  micromagnetic simulations can be numerically expensive, particularly in the case of large-scale finite element simulations. The established  finite element / boundary element method (FEM/BEM) by Fredkin \& Koehler [IEEE Trans.~Mag.~26, 1518 (1990)]  involves a densely populated matrix with unacceptable numerical costs for problems involving a large number of degrees of freedom $N$. By using hierarchical matrices of $\mathcal{H}^2$ type, we show that the memory requirements of this FEM/BEM method can be reduced dramatically, effectively converting the quadratic complexity $\mathcal{O}(N^2)$ of the problem to a linear one $\mathcal{O}(N)$. We obtain matrix size reductions of nearly 99\% in test cases with more than 10$^6$ degrees of freedom, and we test the computed magnetostatic energy values by means of comparison with analytic values. The efficiency of the $\mathcal{H}^2$-matrix compression opens the way to large-scale magnetostatic field calculations in micromagnetic modeling, all while preserving the accuracy of the established FEM/BEM formalism.

\end{abstract}
\end{frontmatter}
\section{Introduction}
Owing to the availability of several public codes with demonstrated reliability, simulating micromagnetic problems based on the Landau-Lifshitz-Gilbert equation \citep{landau_theory_1935,gilbert_phenomenological_2004} has become a standard research task in magnetism. But it remains a challenge to perform large-scale micromagnetic problems involving more than, say, $10^6$ degrees of freedom, especially in the case of non-trivial geometries. In dynamic micromagnetic simulations, where the problem is discretized in time and in space, the effective field contributions of various energy terms \citep{hubert_magnetic_2012} are evaluated frequently, as they need to be updated at each time step. The smooth time integration of the Landau-Lifshitz-Gilbert equation usually requires time steps in the deep sub-ps range, and  simulation studies often address dynamic magnetization processes unfolding over several nano\-seconds. Therefore, in the course of a single micromagnetic simulation, the total number of effective field evaluations may be in the range of hundreds of thousands.
Accordingly, large-scale micromagnetic simulations can only be feasible if they use highly efficient numerical methods to calculate the effective field.

The long-range magnetostatic interaction accounts for the energy due to the dipolar coupling within the ensemble of all microscopic magnetic moments in a ferromagnetic system \cite{hubert_magnetic_2012}. It is usually the most critical effective field contribution in terms of computation time. In its simplest form, which is sometimes referred to as {\em direct integration}, the magnetostatic term results in a twofold volume integral over the magnetostatic charge distribution. In practice, performing such a direct integration is not feasible due to its numerical costs, and corres\-pon\-dingly this approach is almost never used. Instead, one usually aims at solving the Poisson equation for the magnetostatic potential $u$ from which the magnetostatic field results as a gradient field.

While micromagnetic finite difference codes typically employ fast Fourier transforms (FFT) to solve this partial differential equation \citep{berkov_solving_1993}, the irregular mesh used in finite element schemes makes it necessary to use other methods. The Poisson equation in micromagnetics is an open boundary problem \citep{qiushi_chen_review_1997}, meaning that the  boundary conditions are defined by $u$ decaying sufficiently rapidly to zero as the distance from the source approaches infinity. From a technical perspective, a main difficulty lies in the proper consideration of the boundary conditions. Common approaches to treat open boundary conditions in finite element simulations consist in truncating the domain of infinite size at a certain distance, transforming the unbounded space surrounding the sample to a bounded domain by means of bijective spatial projections \citep{brunotte_finite_1992}, or applying a specific boundary integral operator \cite{fredkin_numerical_1990}. In the latter approach, the numerical representation of this integral operator results in a densely populated matrix and thus in a rapid increase of memory requirements with increasing problem size. Hierarchical matrices ($\mathcal{H}$-matrices) \cite{hackbusch_hierarchical_2015} are a powerful method to drastically reduce both, memory requirements and computation time in the finite-element calculation of magnetostatic fields \citep{forster_fast_2003,knittel_compression_2009, kakay_speedup_2010}. Here we show how an improved variant, the $\mathcal{H}^2$-matrix method \citep{hackbusch_khoromskij_sauter_2000,borm_efficient_2010}, can be even more powerful in reducing the computational costs of applying such a discretized integral operator. Owing to this matrix compression, a nearly linear scaling of the problem complexity is achieved. Already at moderate problem sizes, the $\mathcal{H}^2$ compression reduces the memory requirements for the integral operator by about 99\%. In the case of particularly large problems, uncompressed matrices with storage sizes approaching the TB range are thereby reduced to merely a few GB. This dramatic reduction of memory requirements opens the possibility to solve problems of sizes that are unprecedented in micromagnetic finite element simulations. It brings the numerical costs of magnetostatic field calculations in finite elements on a par with those of FFT-based finite-difference schemes, all while preserving the geometric flexibility of the finite element method.

%
%

\section{Basic Equations}
Starting from Maxwell's equation
\be
{\bm\nabla}\cdot{\bm B} = 0
\ee
and the constitutive equation $\bm{B} = \mu_0(\bm{H} + \bm{M})$ relating the solenoidal magnetic induction $\bm{B}$ to the irrotational magnetostatic field $\bm{H}$ and the magnetization $\bm{M}$, one obtains 

\be\label{poiss}
-\Delta u= \rho\quad,
\ee
where $\rho = -{\bm\nabla}\cdot{\bm M}$ and $u$ is the scalar magnetostatic potential with $\bm{H}=-\bm{\nabla}u$. 
The normal derivative of the potential $u$ is discontinuous at the surface $\p\Omega$ of a magnetic material defined in the region $\Omega$:
\be
\frac{\p u^i}{\p\bm{n}} - \frac{\p u^o}{\p{\bm{n}}} = \bm{M}\cdot\bm{n}
\ee
where $o$ and $i$ stand for the outer and inner limit of the derivative at the surface $\p\Omega$, respectively, and $\bm{n}$ is the outward oriented surface normal vector.

To solve for $u$ within the magnetic region $\Omega$, one may make the ansatz to split $u$ into two parts $u = u_1 + u_2$, each with specific properties \citep{fredkin_numerical_1990}. The first part $u_1$ is the solution of the Poisson problem in the 
region $\Omega$:
\be\label{u1Poiss}
\Delta u_1 = \bm{\nabla}\cdot\bm{M}
\ee
with Neumann boundary conditions 
\be\label{u1Neumann}
\frac{\p u_1}{\p\bm{n}} = \bm{M}\cdot\bm{n}
\ee
Outside of $\Omega$ the component $u_1$ is zero, i.e., $u_1(\bm{x}) = 0$ if $\bm{x}\notin\Omega$. The second part, $u_2$, solves the Laplace problem
\be\label{u2diri}
\Delta u_2 = 0 
\ee
The derivative of $u_2$ is continuous along the boundary $\partial\Omega$:
\be 
\frac{\p u_2^i}{\p\bm{n}} - \frac{\p u_2^o}{\p\bm{n}}=0
\ee
The potential $u=u_1+u_2$ is continuous along the boundary, while the individual parts $u_1$ and $u_2$ are generally discontinuous.

This separation of the potential $u$ into two parts allows to solve for $u$ by considering only the region of interest $\Omega$, a task which cannot be achieved directly because the boundary conditions of $u\in\p\Omega$ are unknown. 
The strategy \cite{fredkin_numerical_1990} consists in first solving the Poisson-Neumann problem for $u_1$ in $\Omega$. The solution of the remaining part, $u_2$, requires knowledge of the Dirichlet boundary conditions, i.e., the values of $u_2(\bm{x})$ with $\bm{x}\in\p\Omega$.
These values are obtained via an integral equation that uses the part $u_1$ of the potential. By applying the formalism of hybrid finite element / boundary element formulations \citep{salon_hybrid_1985}, it is straightforward to obtain an integral equation relating the component $u_2$ of the surface potential to the component $u_1$ at the surface. The resulting form is \citep{fredkin_numerical_1990}
 \be\label{mainInt}
u_2(\bm{x}) = \oint u_1(\bm{x'})\frac{\p G(\bm{x},\bm{x}')}{\p\bm{n}}{\rm d}S'
 +\left(\frac{\Psi(\bm{x})}{4\pi}-1\right)u_1(\bm{x})
\ee
with $\bm{x}\in\p\Omega$, where $\Psi(\bm{x})$ is the solid angle subtended at $\bm{x}$. 
In the first term on the right hand side, $G(\bm{x},\bm{x}')$ is Green's function 
\be
G(\bm{x},\bm{x'}) = -\frac{1}{4\pi}\frac{1}{\left\vert\bm{x}-\bm{x'}\right\vert}
\ee
with the property 
\be\label{deltaG}
\Delta G(\bm{x},\bm{x'})=\delta(\bm{x}-\bm{x'})
\ee
where $\delta(\bm{x})$ is the Dirac delta function.
Equation (\ref{mainInt}) provides the Dirichlet boundary conditions needed to solve for $u_2$ in $\Omega$ (eq.~\ref{u2diri}), and thereby to determine $u$.

%
%

\section{Discretization}
In order to solve the system (\ref{mainInt}) numerically, we discretize
it using the collocation approach:
given a surface mesh, we denote the nodes by $(\bm{\xi}_i)_{i=1}^N$ and the
corresponding piecewise linear nodal basis by $(\varphi_i)_{i=1}^N$
with $\varphi_i(\bm{\xi}_j)=\delta_{ij}$.
We are looking for an approximation
\begin{equation*}
  \tilde u_1 = \sum_{j=1}^N u_{1,j} \varphi_j
\end{equation*}
of $u_1$ satisfying the collocation equations
\begin{align*}
  u_2(\bm{\xi}_i)
  &= \oint \tilde u_1(\bm{x'}) \frac{\p G(\bm{\xi}_i,\bm{x}')}{\p\bm{n}}
     \text{d}S'
   + \left(\frac{\Psi(\bm{\xi}_i)}{4\pi}-1\right) \tilde u_1(\bm{\xi}_i)\\
  &= \sum_{j=1}^N m_{ij} u_{1,j}
   = (\underline{\underline{M}} \underline{u_1})_i
\end{align*}
with the matrix $\underline{\underline{M}}\in\bbbr^{N\times N}$ given by
\begin{align}\label{eq:matrix_definition}
  m_{ij} &= \oint \varphi_j(\bm{x'})
           \frac{\p G(\bm{\xi}_i,\bm{x}')}{\p\bm{n}} \text{d}S'
        + \left(\frac{\Psi(\bm{\xi}_i)}{4\pi}-1\right) \delta_{ij}.
\end{align}
If we collect the values of $u_2$ at the collocation points in a
vector $\underline{u}_2\in\bbbr^N$ with
\begin{align*}
  u_{2,i} &= u_2(\bm{\xi}_i),
\end{align*}
we obtain
\begin{align}\label{denseMat}
\underline{u}_2&= \left(\mat{K} + \mat{D}\right)\cdot\underline{u}_1=\mat{M}\cdot\underline{u}_1 
\end{align}
The matrix $\mat{M}$ is densely populated and of size $N\times N$, where $N$ is the number of degrees of freedom at the surface. The diagonal matrix $\mat{D}$ contains the elements in the second term on the right hand side of eq.~(\ref{eq:matrix_definition}).
Analytic forms to calculate the elements of the matrix $\mat{K}$ are reported in Ref.~\citep{lindholm_three-dimensional_1984}.

Solving eqs.~(\ref{u1Poiss}), (\ref{u2diri}) with finite elements  \citep{zienkiewicz_finite_2005} is a numerical standard task once the boundary conditions are known. The discretized representation transforms the problem into a set of linear equations that can be solved using suitable algorithms. If necessary, the process can be accelerated with efficient preconditioning and, in the case of particularly large problems, by employing GPUs to solve the linear system of equations. All matrices involved in this numerical procedure are sparse.

Without matrix compression schemes, the bottleneck in terms of computational costs is the dense matrix $\mat{K}$ in eq.~(\ref{denseMat}), the use of which scales with $N^2$ in terms of both storage requirements and number of arithmetic operations in a matrix-vector product. For large-scale problems, the numerical costs become prohibitively high, especially as far as memory requirements are concerned. Hierarchical matrix methods \citep{hackbusch_hierarchical_2015, noauthor_hierarchical_nodate, noauthor_hlibpro_nodate} have proven to be very efficient in the compression of such matrices resulting from integral operations. They strongly reduce the matrix size and they effectively transform the computational costs from an $\mathcal{O}(N^2)$ scaling to an $\mathcal{O}(N\log N)$ scaling. 

%
%

\section{\texorpdfstring{$\mathcal{H}^2$}{H2}-Matrix Compression}

\begin{figure}
\includegraphics[width=.45\linewidth]{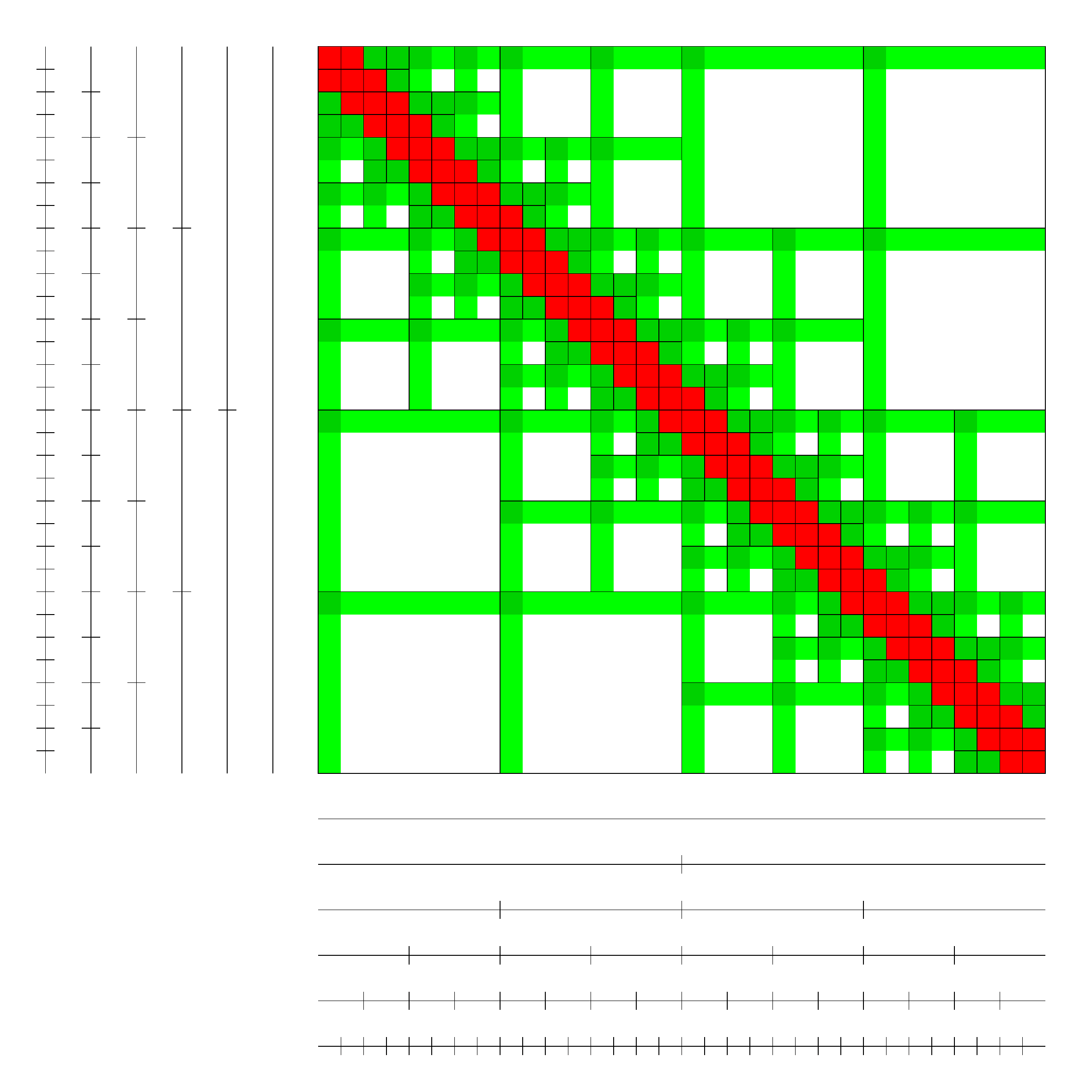}
\qquad
\includegraphics[width=.45\linewidth]{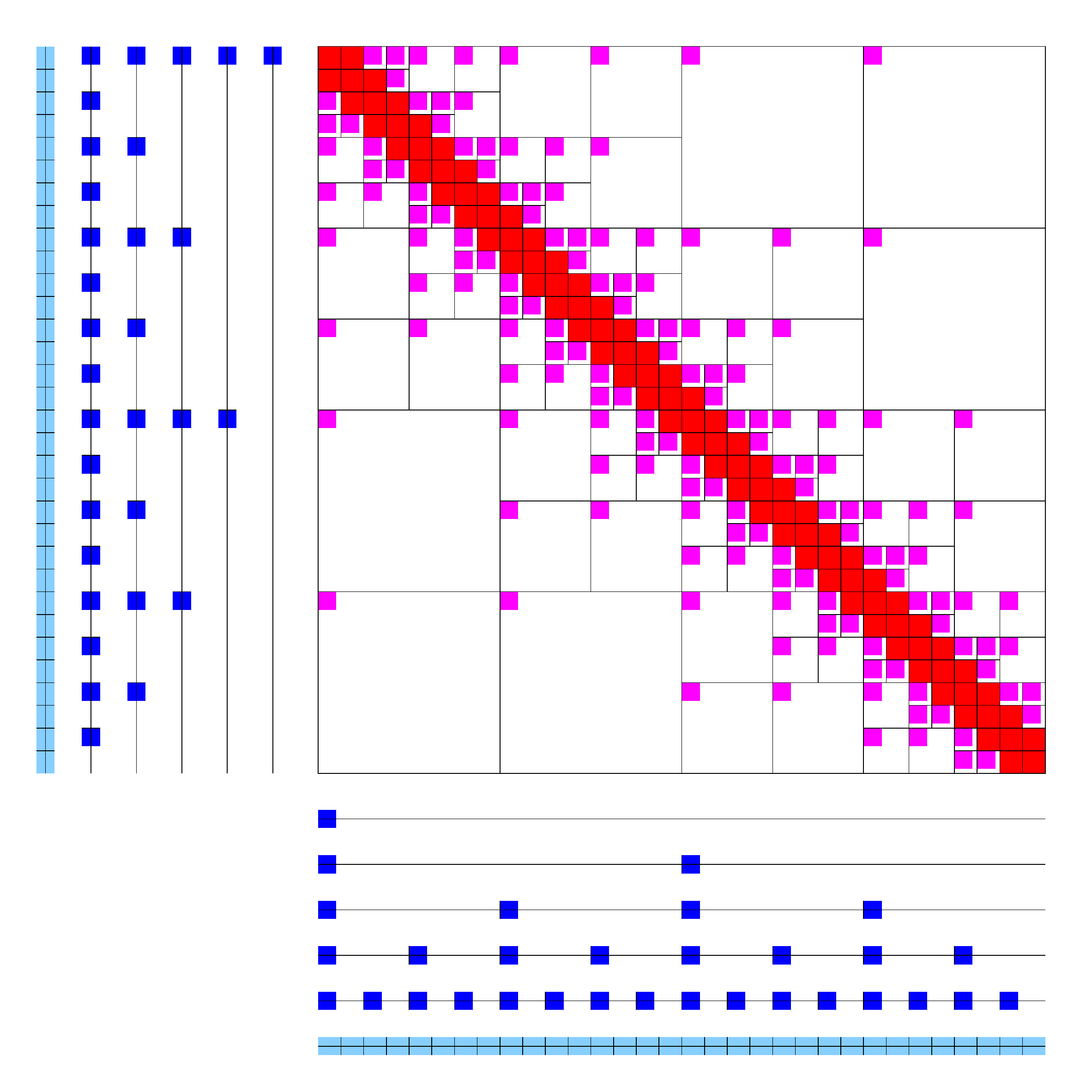}

\caption{\label{hmatrix} Low-rank compression schemes
illustrated for a simple one-dimensional model problem:
$\mathcal{H}$-matrix with $\mathcal{O}(N \log N)$ complexity on the left,
$\mathcal{H}^2$-matrix with $\mathcal{O}(N)$ complexity on the right.}
\end{figure}

Hierarchical matrices are based on the idea that submatrices
$\mat{K}|_{t\times s}$ with index sets $t,s\subseteq[1:N]$
corresponding to well-separated geometrical domains can be approximated
by low-rank matrices and that these low-rank matrices can be represented
in factorized form
\begin{align*}
  \mat{K}|_{t\times s} &\approx \mat{A}\, \mat{B}^T, &
  \mat{A} &\in\bbbr^{t\times k},
  \mat{B} \in\bbbr^{s\times k}
\end{align*}
with the rank $k$.
There are different algorithms for constructing low-rank factorizations
efficiently:
\emph{analytic} methods approximate the kernel function, e.g., by
Taylor expansion, interpolation, or quadrature
\cite{hackbusch_khoromskij_2000,borm_grasedyck_2002,borm_christophersen_2014}
and discretize the result.
These methods are typically very robust and reliable, but they result
in ranks that are higher than necessary.
\emph{Algebraic} methods approximate the matrix $\mat{K}|_{t\times s}$
directly, e.g., by sampling a small number of rows and columns in
the adaptive cross approximation methods \cite{bebendorf_2000}.
These methods are very efficient, but rely on heuristic pivoting
strategies that may lead to unreliable results in certain
situations.
\emph{Hybrid} methods combine an initial analytic approximation with
an algebraic recompression in order to gain the efficiency of
algebraic methods while preserving the reliability of analytic
methods \cite{borm_grasedyck_2004,borm_christophersen_2014}.

In our application, we use the \emph{hybrid cross approximation}
(HCA) \cite{borm_grasedyck_2004} method to construct an $\mathcal{H}$-matrix:
the kernel function $G$ is approximated by Chebyshev interpolation,
the resulting coefficient matrix is compressed using adaptive
cross approximation, and the interpolation is reversed to obtain
the kernel function
\begin{equation*}
  \widetilde{G}(\bm{x},\bm{x}')
  = \sum_{\nu,\mu=1}^k c_{\nu\mu} G(\bm{x},\bm{\eta}'_\mu)
                                  G(\bm{\eta}_\nu,\bm{x}')
\end{equation*}
approximating $G$, where $c_{\nu\mu}$ are suitable scaling coefficients
and $\bm{\eta}_\nu$ and $\bm{\eta}'_\mu$ are interpolation points.
Discretizing this approximation directly yields a low-rank factorization
of $\mat{K}|_{t\times s}$: for $i\in t$ and $j\in s$ with $i\neq j$, we
have
\begin{align*}
  m_{ij} &= \oint \varphi_j(\bm{x}')
           \frac{\p G(\bm{\xi}_i,\bm{x}')}{\p\bm{n}} \text{d}S'\\
  &\approx \sum_{\nu=1}^k \underbrace{\sum_{\mu=1}^k c_{\nu\mu}
           G(\bm{\xi}_i,\bm{\eta}'_\mu)}_{=:a_{i\nu}}
           \underbrace{\oint \varphi_j(\bm{x}')
             \frac{\p G(\bm{\eta}_\nu,\bm{x}')}{\p\bm{n}} \text{d}S'}_{=:b_{j\nu}}.
\end{align*}
Constructing the $\mathcal{H}$-matrix approximation by HCA requires
$\mathcal{O}(N \log N)$ operations and $\mathcal{O}(N \log N)$ units
of storage.
In order to handle very large matrices, we would like to obtain the
optimal linear complexity $\mathcal{O}(N)$.
This is possible by using $\mathcal{H}^2$-matrices
\cite{hackbusch_khoromskij_sauter_2000,borm_hackbusch_2002,borm_efficient_2010},
which use a more sophisticated decomposition that reduces the
storage requirement per block to $k^2$ coefficients.

Using a purely algebraic algorithm \cite[Section~6.5]{borm_efficient_2010},
$\mathcal{O}(N \log N)$ operations are sufficient to ``recompress'' the
$\mathcal{H}$-matrix obtained by HCA into a more efficient
$\mathcal{H}^2$-matrix without a significant loss of accuracy.

Until recently, a wider dissemination of hierarchical matrix methods was impeded by the lack of open-source versions and by licensing restrictions. This has changed with the release of the \texttt{H2Lib} numerical library 
\citep{noauthor_h2lib_nodate} developed in the Scientific Computing Group at Kiel University in Germany. The \texttt{H2Lib} library is an open-source project that provides numerous methods for hierarchical matrices. The modular structure of the library makes it easy to implement these methods in an existing code. In addition to $\mathcal{H}$-matrix algorithms, it also gives the possibility to use $\mathcal{H}^2$-matrices, in particular by offering algorithms to efficiently re-compress $\mathcal{H}$-matrices into $\mathcal{H}^2$-matrices.

We have implemented $\mathcal{H}^2$-matrix compression methods provided by \texttt{H2Lib} to approximate the integral operator of eq.~(\ref{mainInt}) in our micromagnetic finite element code. The use of $\mathcal{H}^2$ methods results in a dramatic reduction of memory requirements, yielding compression factors beyond 99\% and a nearly linear scaling. We demonstrate the superiority of $\mathcal{H}^2$ methods with respect to usual hierarchical matrices in the case of large-scale problems and compare the computed results with analytic values in order to analyze the accuracy.

%
%

\section{Test Systems}
\begin{figure}
\includegraphics[width=.25\linewidth]{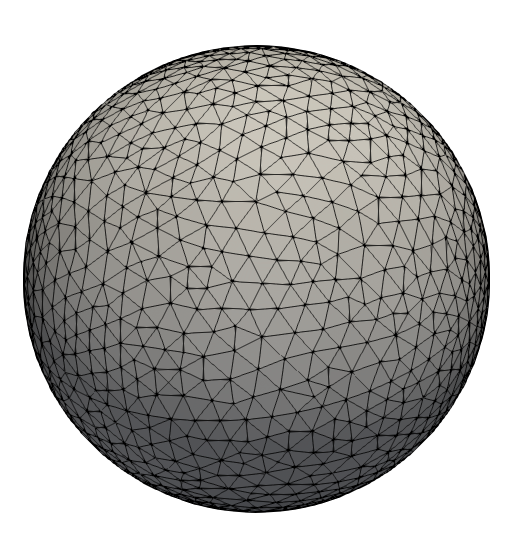}
\includegraphics[width=.4\linewidth]{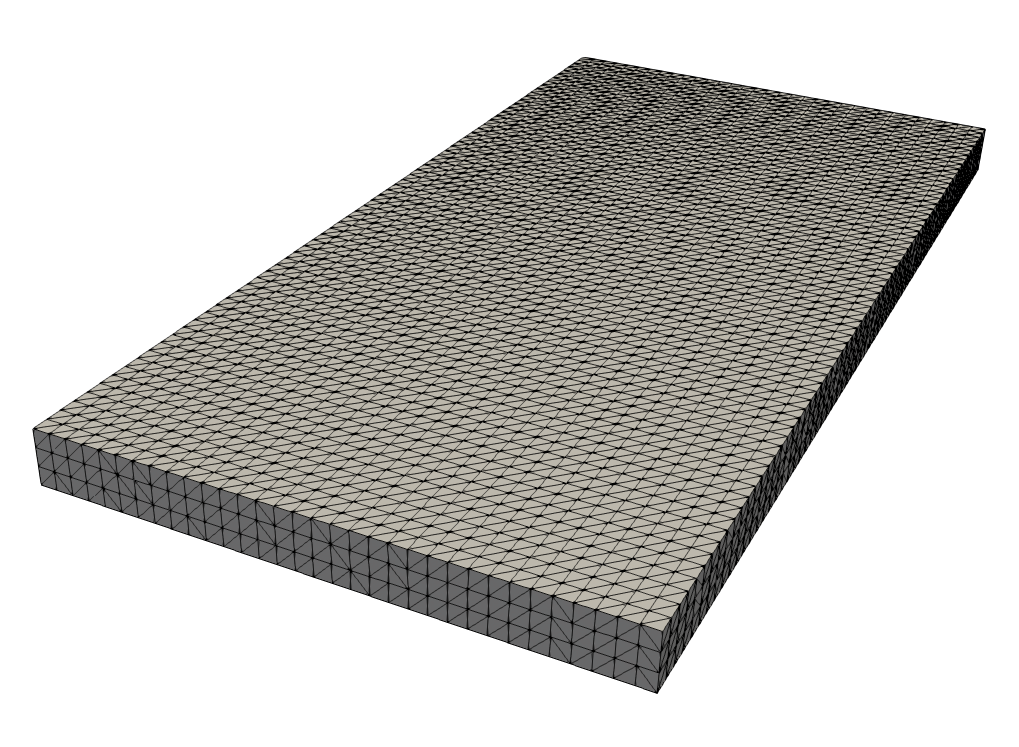}
\includegraphics[width=.33\linewidth]{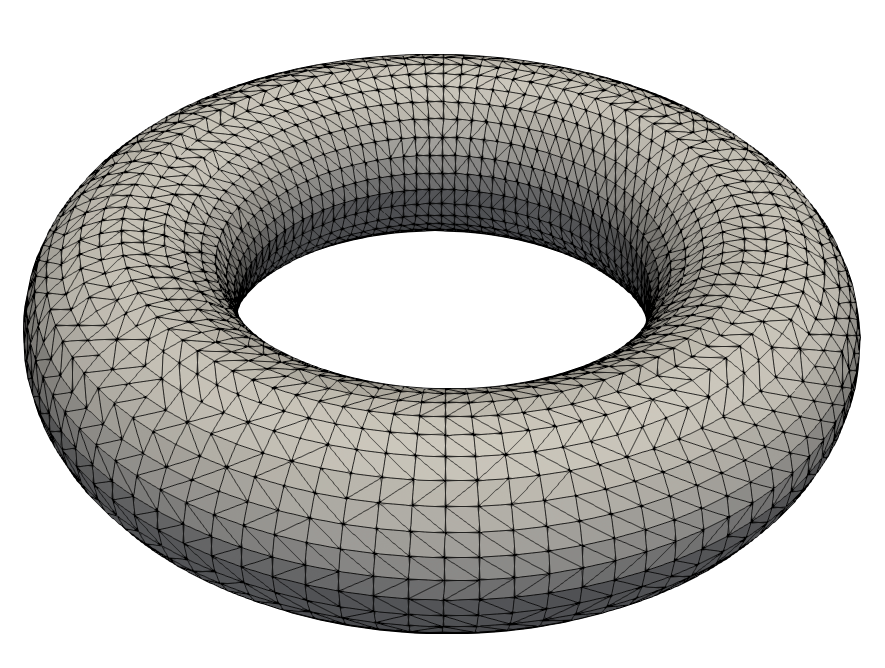}

\caption{\label{geom} Test geometries used to validate the magnetostatic field calculation:  a sphere, a rectangular prism with aspect ratio 10:20:1, and a torus of revolution with aspect ratio of 2.}
\end{figure}

Three different three-dimensional geometries are used to test our implementation of the magnetostatic field calculation with $\mathcal{H}^2$-matrices: a sphere, a rectangular prism, and a torus (Fig.~\ref{geom}). For each of these geometries analytic solutions exist in special cases of the magnetic configuration. This allows us to directly evaluate the accuracy of the computed results. In all cases the absolute size of the problem is irrelevant, as the magnetostatic field calculation is scale invariant.  

%
%

\section{Results}
\subsection{Data compression}
Before discussing the accuracy of the computed results, we first study the effectiveness of the $\mathcal{H}^2$-matrix compression. We compare the matrix sizes as a function of the problem size in three cases: without compression, with $\mathcal{H}$ compression, and with $\mathcal{H}^2$ compression. The data reported here refers to the prism geometry, discretized into linear tetrahedral elements. The other geometries yield very similar data.

In practical micromagnetic simulations, the size of such samples would be in the micron range or above. However, because of the scale invariance of the magnetostatic field calculation, the numerical problem size is unrelated to the physical size. The problem size in the context of matrix compression is mainly determined by the number of surface nodes $N$, since the main task consists in finding a particularly efficient representation of eq.~(\ref{mainInt}), which in its original form relates surface nodes to each other with  $\mathcal{O}(N^2)$ complexity. In contrast to the costs of this surface contribution, the volume part scales in a strictly linear way and is therefore uncritical.

\begin{figure}
\includegraphics[width=\linewidth]{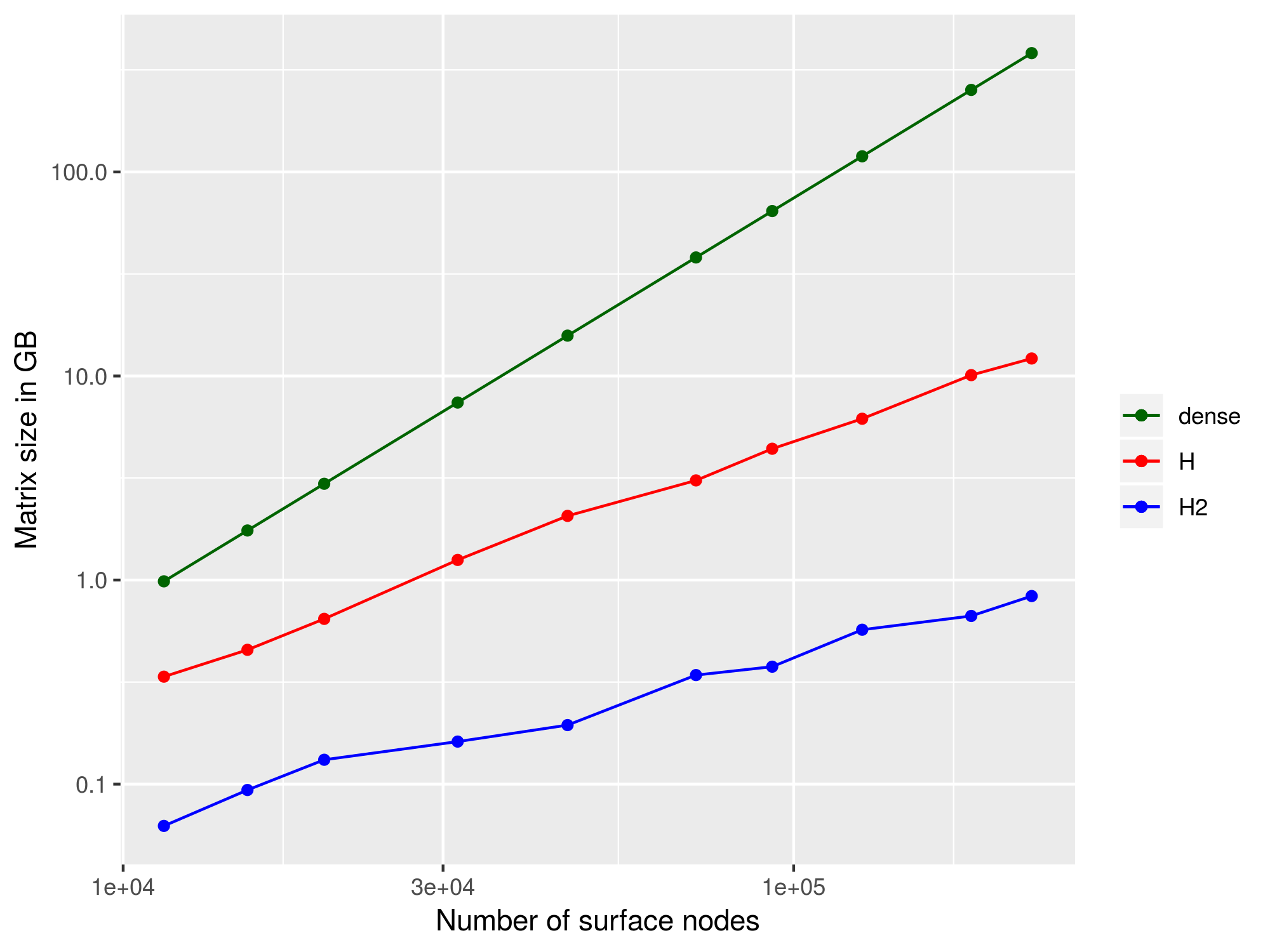}
\caption{\label{matSizes}Comparison of matrix sizes as a function of the number of boundary nodes $N$: uncompressed (dense), $\mathcal{H}$-matrix compression, $\mathcal{H}^2$-matrix compression. Note the double logarithmic scale.}
\end{figure}

In a first step we study how the number of boundary nodes $N$ affects the size of the matrix $\mat{K}$, its compressed version as a traditional hierarchical $\mathcal{H}$ matrix, and its $\mathcal{H}^2$  compressed version. Fig.~\ref{matSizes} displays the results for the prism. The hypothetical size of the uncompressed matrix (green) exceeds 100 GB already at around $N=120\,000$, while the size of the $\mathcal{H}^2$ compressed matrix remains below 1.0 GB even for very large problems with twice as many nodes $N$. 
 The $\mathcal{H}^2$ compression results in matrices that are by a factor of  up to 15 smaller than those obtained with $\mathcal{H}$-matrices.
It should be noted that the number of surface nodes $N$ is typically much smaller than the total number of elements. The data point with the largest number of boundary nodes displayed in Fig.\ref{matSizes} corresponds to a simulation involving more than 13 million finite elements.

\begin{figure}
\includegraphics[width=\linewidth]{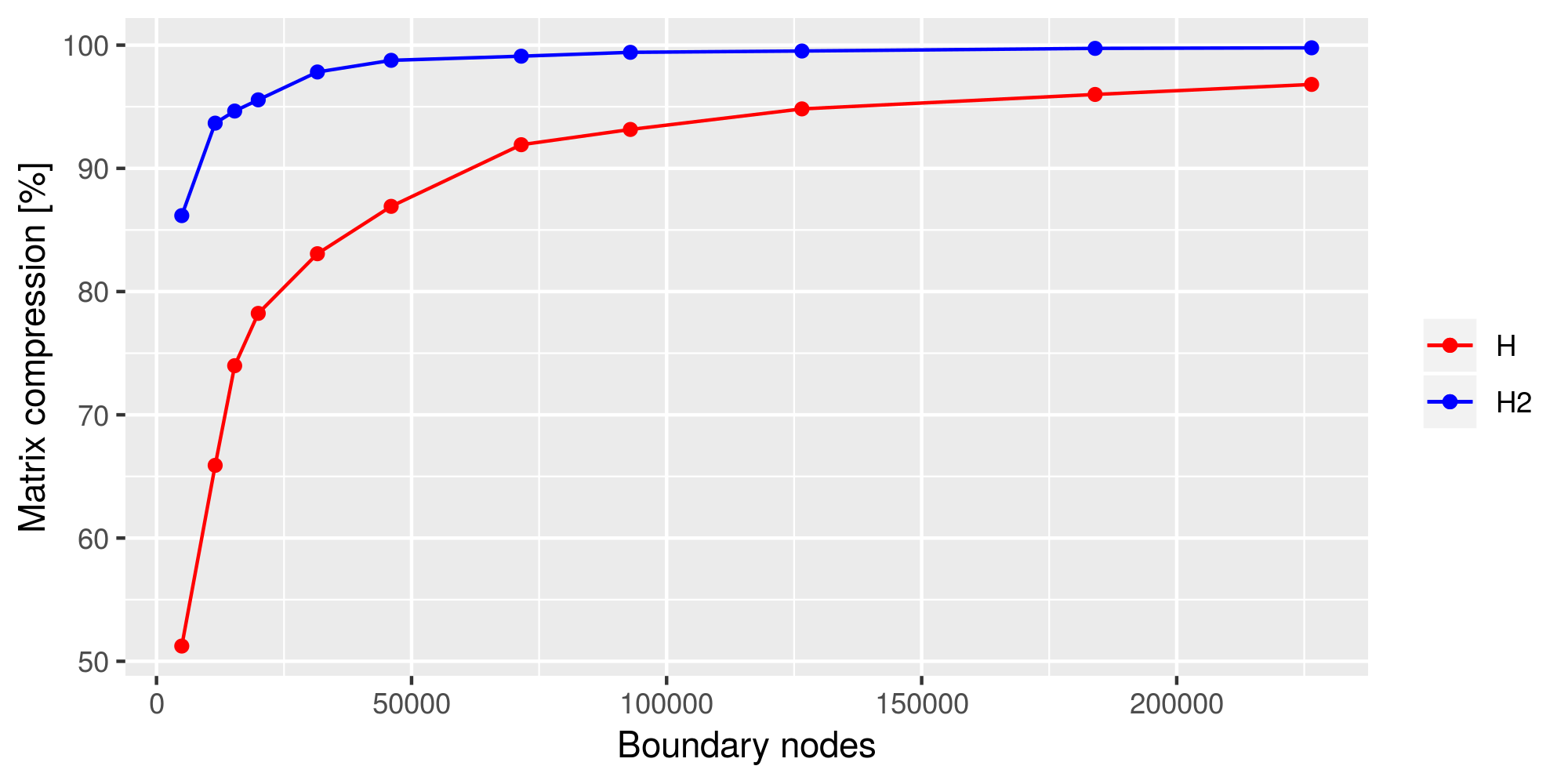}\\
\includegraphics[width=\linewidth]{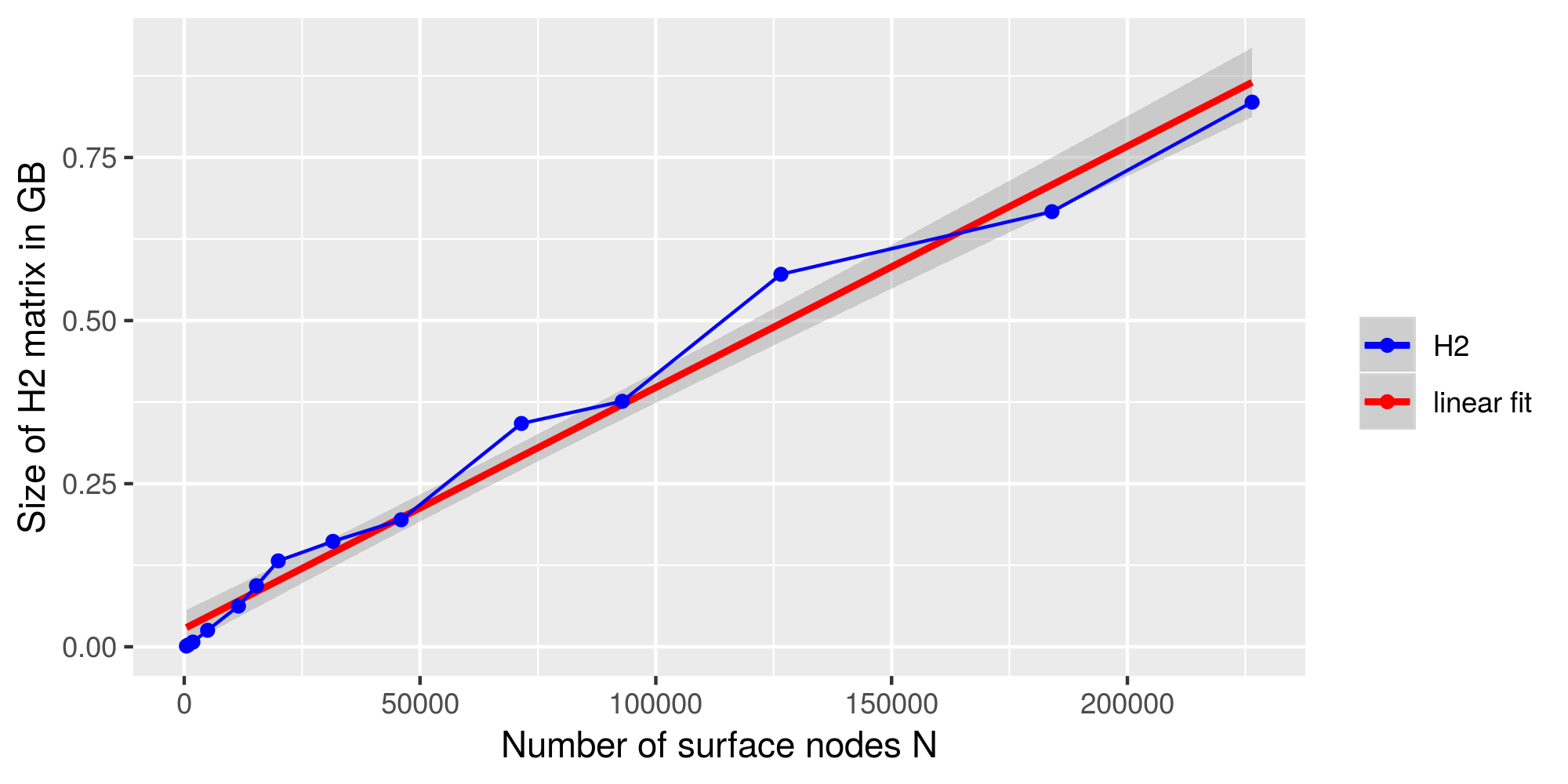}\\
\caption{\label{compressionRatio}Top: matrix compression ratio $r$ (see text) in the case of $\mathcal{H}$- and $\mathcal{H}^2$-matrices. Bottom: the size  of the $\mathcal{H}^2$-matrix scales almost linearly with the number of boundary nodes.}
\end{figure}
Figure \ref{compressionRatio} displays the matrix compression ratio $r$, defined as 
\be
r = 1-\frac{\text{size of compressed matrix}}{
\text{size of dense matrix}}
\ee
for the $\mathcal{H}$- and for the $\mathcal{H}^2$-matrices, highlighting the improved matrix compression provided by the $\mathcal{H}^2$ method.
A virtually linear $\mathcal{O}(N)$ scaling between the size of the $\mathcal{H}^2$-matrix size and the number of boundary nodes $N$ is observed, in contrast to the $\mathcal{O}(N^2)$ scaling of the uncompressed matrix.

In this study we did not analyze the time that is required to set up the $\mathcal{H}^2$-matrices. For our purposes, the speed with which the $\mathcal{H}^2$-matrices are being generated is not important, as this calculation is required only once during the preprocessing. Once the matrix is calculated, it is stored in a file that is read at the beginning of any subsequent calculation for a given geometry. Even for the largest problems that we have analyzed, setting up the $\mathcal{H}^2$-matrix did not exceed the range of several minutes of computation. The \texttt{H2Lib} package provides optimized methods to speed up the calculation of these matrices even further by directly setting up the $\mathcal{H}^2$ matrix, without calculating and compressing an $\mathcal{H}$ matrix.

\subsection{Magnetostatic field and energy}
We now compare the computed values of the magnetostatic energy with analytical ones. Such analytic solutions exist in micromagnetic theory for a few specific geometries with certain magnetic configurations. Homogeneously magnetized ellipsoids are one such example \cite{osborn_demagnetizing_1945}, with the sphere being the simplest case. 

Another case where analytic solutions exist are ferromagnetic rectangular prism, with homogeneous magnetization along any of the main axes. This is a hypothetical situation, as the equilibrium structure of the magnetization can never be homogeneous in non-ellipsoidal particles. Accordingly, we do not relax the magnetization in the calculations, as would be usually done by means of a time integration of the Landau-Lifshitz-Gilbert equation. The resulting, converged state would be useless for our purposes because it would not be possible to compare the result with an analytic solution. 

A trivial analytic solution for the magnetostatic field and energy exists in the case of a torus with purely azimuthal magnetic configuration. In this case, the divergence of the magnetic structure is exactly zero, thereby avoiding magnetostatic volume charges. There are also no surface charges, as the magnetization is perpendicular to the surface normal vector at any point of the surface. This special situation, which is only possible in a multiply connected geometry of revolution, removes any source of the magnetostatic field, so that the expected magnetostatic energy density is exactly zero, $e_d=0$. Also in this case, the purely azimuthal magnetization structure is imposed as an initial condition of the calculation, and it is not relaxed in order to preserve the magnetization state with exact analytic solution.

We calculate the magnetostatic field and determine the magnetostatic field energy for these three special cases with the method described before, which uses a $\mathcal{H}^2$-matrix representation of the integral operator in eq.~(\ref{mainInt}). The Poisson and Laplace equations (\ref{u1Neumann}), (\ref{u2diri}) are solved with a standard sparse LU-solver \cite{eigenweb} for smaller problem sizes, while for larger problems (typically involving more than $10^6$ elements) we used a GPU-accelerated solver using the biconjugate gradient stabilized method (BiCGStab) \cite{cusp} .

\subsubsection{Rectangular prism}
In the case of a rectangular prism, demagnetizing factors can be defined for the case of homogeneous magnetization along the principal axes \cite{aharoni_demagnetizing_1998}. Although the demagnetizing field of a homogeneously magnetized prism is not homogeneous (it displays, on the contrary, pronounced inhomogeneities at the edges and corners), these so-called ballistic demagnetizing factors can be used to calculate the magnetostatic energy density. For example, in the case of homogeneous magnetization along the $x$ axis (assuming that it is one of the prism axes), analytic expressions for $N_x(a,b,c)$ exist as a function of the prism's aspect ratio $a:b:c$ \cite{aharoni_demagnetizing_1998}. These demagnetizing factors allow us to calculate the spatially averaged energy density $e_d= \mu_0N_xM_s^2/2$.
In our case we consider a rectangular prism with aspect ratio of 10:20:1 with homogeneous out-of-plane magnetization. The computed magnetostatic energy can be compared with the analytic value derived from the demagnetizing factors of a rectangular prism with these axis ratios.

\begin{figure}
\includegraphics[width=\linewidth]{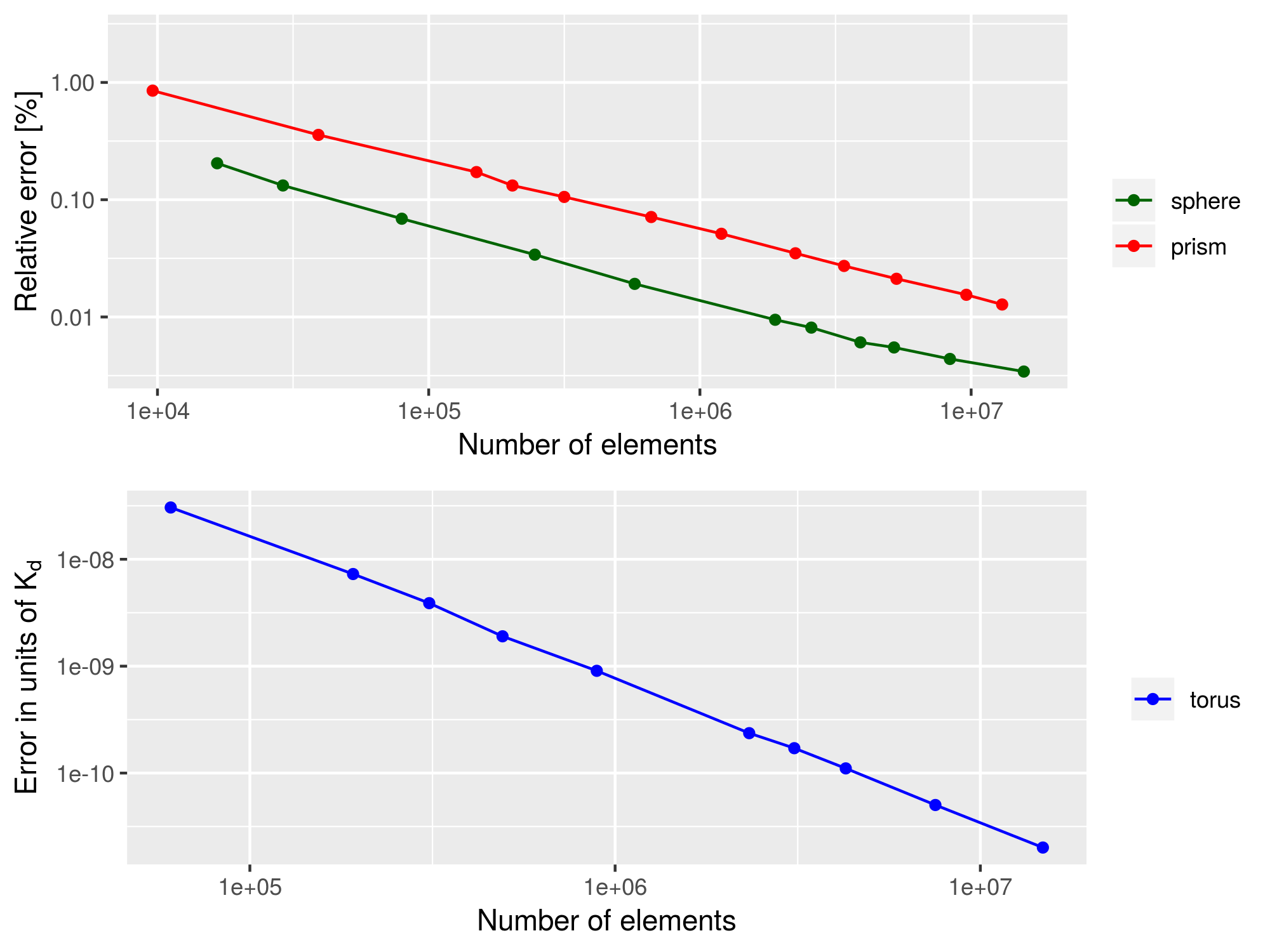}
\caption{\label{err_geoms}Top: Relative error of the calculated magnetostatic energy of a homogeneously magnetized sphere (green) and a prism (red) as a function of the number of elements. Bottom: Absolute error of the magnetostatic energy for a torus of revolution with aspect ratio of 2 and azimuthal magnetization. The error is displayed in units of the stray field constant $K_d=\mu_0M_s^2/2$.}
\end{figure}

\subsubsection{Sphere}
For a homogeneously magnetized sphere, the demagnetizing field is homogeneous, so that it can be related to the magnetization field $\bm{H}_d$ simply by means of multiplication with a scalar demagnetizing factor $\bm{H}_d = -N \bm{M}$, where $N = 1/3$. The magnetostatic energy density $e_d$ of a homogeneously magnetized sphere with saturation magnetization $M_s$ is $e_d = \mu_0M_s^2$/6. The energy resulting from the simulation should match this analytic value, with any deviation from it representing a quantifiable error.

\subsubsection{Torus}
In the case of a torus, the expected magnetostatic energy density $e_d$ is zero, as discussed before. It is therefore not possible to calculate the relative error of the numerically computed field. We display instead the unsigned absolute error in Fig. ~\ref{err_geoms} for a magnetic torus with aspect ratio 2:1 and purely azimuthal magnetization. The data is represented here in a dimensionless way, in units of the stray field constant \cite{hubert_magnetic_2012} defined as $K_d=\mu_0M_s^2/2$. This normalization makes it possible to put the data in a context of general validity, such that it is not specific to the material type used in our simulations. 

\subsubsection{Numerical errors}
In all geometries considered here, the corresponding analytic value is nicely reproduced, and the numerical error decreases rapidly with increasing discretization density (Fig.~\ref{err_geoms}). For the rectangular prism we observe a lower convergence rate and a somewhat larger error compared to the other geometries. 
This is due to the different nature of the solutions of the demagnetizing field. While the demagnetizing field of a homogeneously magnetized sphere is homogeneous, a rectangular prism with homogeneous magnetization results in a strongly inhomogeneous $\bm{H}$ field displaying a logarithmic divergence at the edges and corners \cite{rave_corners_1998,thiaville_corner_1998}. With a piecewise linear finite element discretization, it is much more challenging to obtain a good approximation of such a strongly inhomogeneous field than it is to model a homogeneous field, as it result in an ellipsoid. For the torus geometry with azimuthal magnetization the solution is even simpler to approximate than the field of a homogeneously magnetized sphere, es in that case the anlaytic result is a vanishing field $\bm{H}$. 

All cases indicate that the global numerical error is not due to the approximation connected with the $\mathcal{H}^2$ matrix compression. In the case of the prism the error is most likely due to inaccuracies in the discretized representation of the solution of the Poisson and Laplace equation (eqs.~\ref{u1Poiss}, \ref{u2diri}) within the volume, primarily because of pronounced inhomogeneities at the corners, while in the case of the sphere and the torus the quality of the geometric approximation of the curved surfaces by means of piecewise linear segments is likely to be a main source of errors. In any case, it should be noted that the error remains very small, perfectly in the range of acceptable limits for micromagnetic simulations.
Also, the logarithmic divergences in the case of the prism result from an artificial situation where a homogeneous magnetization is assumed. These divergences are removed in realistic magnetic structures, where specific magnetic configurations develop at the edges and corners \cite{rave_magnetic_1998}.


%
%

\section{Conclusions}
We have implemented the $\mathcal{H}^2$-matrix formalism to solve the magnetostatic open-boundary problem in large-scale micromagnetic simulations. The implementation of the {\tt H2Lib} library with its powerful methods to compress dense matrices that occur in numerical forms of integral operators has resulted in compression ratios of nearly 99\%, matrix sizes that are up to a factor of 15 smaller compared to traditional hierarchical matrices, and a nearly linear $\mathcal{O}(N)$ scaling of the memory requirements. By combining the $\mathcal{H}^2$ method with GPU-accelerated solvers we have simulated the magnetostatic field in three basic geometries with very large discretization density and compared the result to analytic solutions.
The results demonstrate that the method is accurate and efficient for large-scale problems. 

We have discussed the advantages of using the $\mathcal{H}^2$-method in micromagnetic simulations in terms of matrix compression and accuracy, but not in terms of computation time. Although we can ascertain that the method is very fast, in any case much faster than the use of dense matrices, a detailed benchmarking of the calculation speed of the FEM/BEM method with $\mathcal{H}^2$-matrices is beyond the scope of this article. A precise analysis of the speed of the algorithm is a complicated issue because the scheme involves several different steps, each of which could be subject to optimization. The number of arithmetic operations involved in a matrix-vector multiplication with $\mathcal{H}^2$-matrices shows nearly linear complexity, as does the computation time. In the hybrid FEM-BEM method used here, the $\mathcal{H}^2$-matrix-vector multiplication is only one part of a larger procedure, and this part is not critical for the overall computation speed. If efficient matrix compression schemes like those offered by {\tt H2Lib} are used, it is the iterative solution of linear equations that takes  the role of the bottleneck in terms of computation time in large problems. This latter part can be greatly accelerated with suitable preconditioning and by means of GPU parallelization. 

In conclusion, the implementation of efficient matrix compression methods and the resulting reduction of memory requirements provides an important ingredient in the development of new, powerful micromagnetic finite-element codes. It makes the  simulation of realistic situations accessible which currently cannot be solved with the available standard codes, e.g., cases involving large sample sizes with non-trivial geometry or large arrays of interacting magnetic particles.

%
%

\section*{References}
\bibliography{h2bib}

\end{document}